\documentclass[reqno]{amsart}
\usepackage{setspace}
\usepackage{graphicx}
\usepackage{amsmath}
\usepackage{amsfonts}
\usepackage{caption}
\usepackage{amssymb}
\usepackage{longtable}
\usepackage[margin=2.5cm]{geometry}
\usepackage{float}

\begin{document}

\centerline{\Huge \bf New Lattices with Unexpected}
\centerline{\Huge \bf Shallow Deep Holes}

\bigskip
\centerline{\Large Yiming Li$^a$ and Chuanming Zong$^{a,b,c}$}

\bigskip
\centerline{$^a$Center for Applied Mathematics, Tianjin University, Tianjin 300072, China}

\centerline{$^b$School of Cybersecurity, Tianjin University, Tianjin 300072, China}

\centerline{$^c$Institute of Post-Quantum Science and Technology, Tianjin University of Technology}
\centerline{ Tianjin 300384, China}

\medskip
\centerline{Correspondence: xiaozhuang@tju.edu.cn; cmzong@tju.edu.cn}

{\large
\vspace{0.8cm}
\centerline{\begin{minipage}{13.4cm}
{\bf Abstract.} Among lattices normalized to have packing radius one, the known shallowest depths of the deep holes in $7$-, $8$-, $9$- and $10$-dimensional lattices are $\sqrt{7/3}$, $\sqrt2$, $\sqrt{5/2}$ and $\sqrt{8/3}$, given by $E_7^*$, $E_8$, the laminated lattices $\Lambda_9$ and $\Lambda_{10}$, respectively. This paper presents new lattices in these dimensions with shallower deep holes. It also answers two problems of J. Martinet about short representatives of the classes of $\Lambda/2\Lambda$ in $7,8,9,10$ and $12$ dimensions.
\end{minipage}}

\bigskip\noindent
2020 {\it Mathematics Subject Classification}: 11H06, 11H31, 11H50, 52C17.

\medskip
\noindent
{\bf Keywords:} lattice; deep hole; packing-covering constant.

\vspace{0.8cm}
\noindent
{\LARGE\bf 1. Introduction}

\bigskip\noindent
Let $\{{\bf a}_1,{\bf a}_2,\ldots,{\bf a}_n\}$ be $n$ linearly independent vectors in the $n$-dimensional Euclidean space $\mathbb{E}^n$. We call
$$\Lambda=\{z_1{\bf a}_1+z_2{\bf a}_2+\ldots+z_n{\bf a}_n:\ z_i\in\mathbb{Z}\}$$
an $n$-dimensional {\it lattice} and call $\{{\bf a}_1,{\bf a}_2,\ldots,{\bf a}_n\}$ a {\it basis} of the lattice $\Lambda$. Lattice is a basic concept in mathematics, which has been investigated by Gauss, Hermite, Dirichlet, Minkowski, Voronoi, Korkin, Zolotarev and many modern authors from various respects (see \cite{Co01, Zo99}). Recent years, it has become the key for post-quantum cryptography (see \cite{Zo25}).

Let $B^n$ denote the unit ball centered at the origin ${\bf o}$ of $\mathbb{E}^n$, let $R(\Lambda)$ denote the smallest number $R$ such that $RB^n+\Lambda$ is a covering of $\mathbb{E}^n$, and let $r(\Lambda)$ denote the largest number $r$ such that $rB^n+\Lambda$ is a packing in $\mathbb{E}^n$. Usually, the points in $\mathbb{E}^n$ with exact distance $R(\Lambda)$ to $\Lambda$ are called the {\it deep holes} of $\Lambda$ (see \cite[p.33]{Co01}).  In 1950, C. A. Rogers \cite{Ro01} defined
$$\gamma^*(B^n)=\min_{\Lambda\in\mathcal{L}_n}\frac{R(\Lambda)}{r(\Lambda)},$$
where $\mathcal{L}_n$ denotes the family of all $n$-dimensional lattices. Usually, it is called the {\it packing-covering constant}. By an inductive method, Rogers proved that
$$\gamma^*(B^n)\leq 3$$
holds in every dimension. In 1972, via mean value techniques developed by Rogers \cite{Ro02} and C. L. Siegel \cite{Si01}, G. L. Butler \cite{Bu01} improved Rogers' upper bound to
$$\gamma^*(B^n)\leq 2+o(1).$$
In fact, these upper bounds hold for all $n$-dimensional centrally symmetric convex bodies (see \cite{Zo02}).

Clearly, $2r(\Lambda )$ is the length of the shortest nonzero vectors of $\Lambda$. It is well-known that to determine the shortest vectors and the smallest covering radius $R(\Lambda)$ of a given lattice $\Lambda$ are two fundamental problems in lattice-based cryptography (see \cite{Zo26}). Naturally, Rogers' constants $\gamma^*(B^n)$ have applications in lattice-based cryptography (see \cite{Mi04}).

As one can imagine, similar to the packing density problem and the covering density problem, to determine the exact value of $\gamma^*(B^n)$ is a challenging job. The known exact results and the known best upper bounds are listed in the following tables:

\begin{table}[ht]
\newcommand{\AD}[1]{{\begin{tabular}{@{}c@{}}#1\end{tabular}}}
\begin{tabular}{|c|c|c|c|c|c|}
\hline $n$ & $2$ & $3$ & $4$ & $5$ \\
\hline $\gamma^*(B^n)$ & $2/\sqrt3$ & $\sqrt{5/3}$ & $\sqrt{2\sqrt3}(\sqrt3-1)$ & $\sqrt{3/2+\sqrt{13}/6}$\\
\hline Authors & Trivial & K. B\"{o}r\"{o}czky \cite{Bo01} & J. Horv\'{a}th \cite{Ho01} & J. Horv\'{a}th \cite{Ho02}\\
\hline
\end{tabular}

\vspace{0.3cm}
\centerline{Table 1. The known exact results for $\gamma^*(B^n)$.}
\end{table}

\begin{table}[ht]
\newcommand{\AD}[1]{{\begin{tabular}{@{}c@{}}#1\end{tabular}}}
\begin{tabular}{|c|c|c|c|c|c|c|c|}
\hline $n$ & $6$ & $7$ & $8$ & $9$ & $10$ & $24$ \\
\hline $\gamma^*(B^n)$ & $\leq2\sqrt{2\sqrt{798}-56}$ & $\leq\sqrt{7/3}$ & $\leq\sqrt2$ & $\leq\sqrt{5/2}$ & $\leq\sqrt{8/3}$ & $\leq\sqrt2$\\
\hline
\end{tabular}

\vspace{0.3cm}
\centerline{Table 2. The known best upper bounds for $\gamma^*(B^n)$.}
\end{table}

\noindent
The upper bounds in Table 2 are given by $\Lambda(Q^{pc}_6)$, $E^*_7$, $E_8$, the laminated lattices $\Lambda_9$, $\Lambda_{10}$ and the Leech lattice $\Lambda_{24}$, respectively, see \cite{Co02, Co01, Sc02, Sc01, Wo01}. In 2002, C. Zong \cite{Zo02} conjectured that
$$\gamma^*(B^8)=\sqrt2,\qquad \gamma^*(B^{24})=\sqrt2,$$
where equalities hold at $E_8$ and $\Lambda_{24}$, respectively.

\medskip
This paper mainly proves the following theorem.

\medskip
\noindent
{\bf Theorem 1.} {\it In $7$-$10$ dimensions, the upper bounds in Table $2$ can be improved as
$$\gamma^*(B^7)\leq1.51582801\cdots<\sqrt{7/3},$$
$$\gamma^*(B^8)\leq1.41413038\cdots<\sqrt2,$$
$$\gamma^*(B^9)\leq1.57982271\cdots<\sqrt{5/2}$$
and
$$\gamma^*(B^{10})\leq1.56347674\cdots<\sqrt{8/3}.$$}

In fact, we prove Theorem 1 by constructing lattices whose short vectors determine their Voronoi cells and then computing their covering radii. These constructions also provide affirmative answers to the following problems proposed by J. Martinet in 2002 in corresponding dimensions.

\medskip
\noindent
{\bf Problem 1 (Martinet \cite{Ma01}).} {\it Do there exist lattices in dimensions $n\geq7$ such that all classes in $\Lambda/2\Lambda$ possess representatives of norm $N<8r(\Lambda)^2$?}

\medskip
\noindent
{\bf Problem 2 (Martinet \cite{Ma02}).} {\it Do there exist lattices in dimensions $n>10$, not similar to $K_{12}$ (Coxeter-Todd lattice) nor to $\Lambda_{24}$ (Leech lattice), such that all classes in $\Lambda/2\Lambda$ possess representatives of norm $N\leq 8r(\Lambda)^2$?}

\medskip
More precisely, to approach Martinet's problems, this paper proves the following theorems.

\medskip
\noindent
{\bf Theorem 2.} {\it For $n=7, 8, 9$ and $10$, there exists an $n$-dimensional lattice $\Lambda$ such that all classes in $\Lambda/2\Lambda$ possess representatives of norm $N<8r(\Lambda)^2$.}

\medskip
\noindent
{\bf Theorem 3.} {\it There exists a $12$-dimensional lattice $\Lambda$ (not similar to $K_{12}$) such that all classes in $\Lambda/2\Lambda$ possess representatives of norm $N\leq 8r(\Lambda)^2$.}

\medskip
\noindent
{\bf Remark 1.} It is well-known that $E_6$, $E_7$ and $E_8$ provide the densest lattice sphere packing in corresponding dimensions (see \cite{Co01, Zo99}). On the other hand, it has been shown by A. Sch\"{u}rmann and F. Vallentin \cite{Sc02, Sc01} that these lattices do not provide the thinnest lattice sphere covering in corresponding dimensions.

\medskip
\noindent
{\bf Remark 2.} In this paper, we often use numerical computation. Whenever we claim a routine computation by computer we mean simple numerical computation based on Matlab.

\vspace{0.8cm}
\noindent
{\LARGE\bf 2. Preliminaries}

\bigskip
\noindent
For an $n$-dimensional lattice $\Lambda$, define
$$X(\Lambda)=\big\{{\bf v}\in\Lambda:\ 2r(\Lambda)\leq||{\bf v}||<2\sqrt2r(\Lambda)\big\},$$
where $||{\bf v}||$ denotes the length of vector ${\bf v}$.

\medskip
\noindent
{\bf Lemma 1 (Li and Zong \cite{LiY01}, Martinet \cite[p.867]{Ma01}).} {\it For an $n$-dimensional lattice $\Lambda$, we have}
$${\rm card}\{X(\Lambda)\}\leq 2^{n+1}-2.$$

For an arbitrary vector ${\bf v}\in\mathbb{E}^n\setminus\{{\bf o}\}$, define the half space
$$H_{{\bf v}}=\big\{{\bf x}\in\mathbb{E}^n: ||{\bf x}||\leq||{\bf x}-{\bf v}||\big\}.$$
The {\it Voronoi cell} of an $n$-dimensional lattice $\Lambda$ is defined as the set
$$V(\Lambda)=\cap_{{\bf v}\in\Lambda\setminus\{{\bf o}\}}H_{{\bf v}}.$$
It is easy to see that $V(\Lambda)$ is a centrally symmetric convex polytope and can tile the space $\mathbb{E}^n$ by lattice $\Lambda$. Therefore, we have
$$R(\Lambda)=\max\{||{\bf x}||:\ {\bf x}\ {\rm is\ a\ vertex\ of}\ V(\Lambda)\}.$$

We call the minimal set of lattice vectors $M(\Lambda)\subset\Lambda$ such that
$$V(\Lambda)=\cap_{{\bf v}\in M(\Lambda)}H_{{\bf v}}$$
the set of {\it Voronoi relevant} vectors.

\medskip
\noindent
{\bf Lemma 2 (Voronoi (see \cite[p.477]{Co01})).} {\it Let $\Lambda$ be a lattice and ${\bf v}\in\Lambda\setminus\{{\bf o}\}$ a nonzero lattice vector. Then ${\bf v}$ is Voronoi relevant if and only if $\pm{\bf v}$ are the only two shortest vectors in the coset $2\Lambda+{\bf v}$.}

\medskip
\noindent
{\bf Remark 3.} If an $n$-dimensional lattice $\Lambda$ satisfies
$${\rm card}\{X(\Lambda)\}=2^{n+1}-2,$$
then the equality implies that every nonzero class of $\Lambda/2\Lambda$ contains exactly one pair $\pm {\bf v}$ from $X(\Lambda)$ (see \cite{LiY01, Ma01}). This pair consists of the only shortest vectors in that class, so Lemma 2 shows that the set $M(\Lambda)$ of Voronoi relevant vectors of the lattice $\Lambda$ satisfies
$$M(\Lambda)=X(\Lambda).$$

Given a lattice covering of a convex body $C$ in $\mathbb{E}^n$, its star number is the number of the translates of $C$ by lattice vectors, including $C$, which intersect the body $C$. In 1964, P. Erd\"{o}s and Rogers \cite{Er01} proved the following result.

\medskip
\noindent
{\bf Lemma 3.} {\it Let $C$ be an ${\bf o}$-symmetric strictly convex body and $\Lambda$ a covering lattice of $C$ in $\mathbb{E}^n$. Then the star number of the covering $\{C+{\bf v}:\ {\bf v}\in  \Lambda\}$ is at least $2^{n+1}-1$.}

\medskip
Combining Lemma 1 and Lemma 3, Y. Li and C. Zong \cite{LiY01} obtained the following necessary condition.

\medskip
\noindent
{\bf Lemma 4.} {\it If an $n$-dimensional lattice $\Lambda$ satisfies
$$R(\Lambda)/r(\Lambda)<\sqrt2,$$
then every class in $\Lambda/2\Lambda$ contains a representative of norm
$$N<8r(\Lambda)^2.$$
In other words, if no lattice satisfying the condition in Problem 1 exists in a dimension $m$, then we have}
$$\gamma^*(B^m)\geq\sqrt2.$$

By Lemma 4, we can see that Martinet's problems provide a useful starting point for studying the deep hole problems, particularly in dimensions 8 and 24.

It is well known that there is a correspondence between lattices and positive definite quadratic forms (see \cite[p.23]{Zo99}). On the one hand, for each $n$-dimensional lattice
$$\Lambda=A\mathbb{Z}^n,$$
where $A$ is a nonsingular $n\times n$ matrix, there is a corresponding positive definite matrix
$$G=A^{\mathsf T}A$$
and the associated quadratic form
$$Q({\bf z})={\bf z}^{\mathsf T}G{\bf z}.$$
On the other hand, for each positive definite matrix $G$ and the associated quadratic form $Q({\bf z})={\bf z}^{\mathsf T}G{\bf z}$, there is a corresponding lattice
$$\Lambda=A\mathbb{Z}^n$$
such that
$$G=A^{\mathsf T}A.$$

We call $G$ the {\it Gram matrix} of the basis formed by the columns of $A$. Clearly, for an arbitrary lattice point
$${\bf v}=A{\bf z}\in\Lambda,$$
the norm of ${\bf v}$ is
$$||{\bf v}||^2=Q({\bf z}).$$
Therefore, we have
$$r(\Lambda)^2=\frac{1}{4}\min_{{\bf z}\in\mathbb{Z}^n\setminus\{{\bf o}\}}Q({\bf z}).$$

\medskip
\noindent
{\bf Remark 4.} For an arbitrary vector ${\bf x}=A{\bf y}$, we have
$$||{\bf x}||\leq||{\bf x}-A{\bf z}||$$
if and only if
$$2{\bf y}^{\mathsf T}G{\bf z}\leq{\bf z}^{\mathsf T}G{\bf z}.$$
Therefore, if an $n$-dimensional lattice $\Lambda$ satisfies
$${\rm card}\{X(\Lambda)\}=2^{n+1}-2,$$
by Remark 3 we have
$$V(\Lambda)=AP^*$$
and
$$R(\Lambda)^2=\max\big\{{\bf y}^{\mathsf T}G{\bf y}:\ {\bf y}\ {\rm is\ a\ vertex\ of\ } P^*\big\},$$
where
$$P^*=\big\{{\bf y}\in\mathbb{E}^n:\ 2{\bf y}^{\mathsf T}G{\bf z}\leq{\bf z}^{\mathsf T}G{\bf z}\ {\rm for\ all\ {\bf z}\in \mathbb{Z}^n\ satisfying\ } 4r(\Lambda)^2\leq{\bf z}^{\mathsf T}G{\bf z}<8r(\Lambda)^2\big\}.$$

If $G$ has rational entries, the vertices of $P^*$ can be enumerated exactly using the reverse search
algorithm described in \cite{Av01}. For each vertex ${\bf y}$ of $P^*$, the corresponding Voronoi vertex $A{\bf y}$ has norm ${\bf y}^{\mathsf T}G{\bf y}$. Consequently, $R(\Lambda)^2$ can be evaluated using exact rational arithmetic.

\medskip
In section 3 and section 4, we present a collection of positive definite quadratic forms in $7,8,9,10$ and $12$ dimensions, which proves Theorem 1, Theorem 2 and Theorem 3.

\vspace{0.8cm}
\noindent
{\LARGE\bf 3. Proof of Theorem 1 and Theorem 2}

\bigskip
\noindent
{\Large\bf 3.1. The Seven-Dimensional Case}

\medskip
Let $Q_7({\bf z})={\bf z}^{\mathsf T}G_7{\bf z}$, where $G_7=N_7/d_7,\ d_7=999999999739$ and
            \[
            \resizebox{\textwidth}{!}{%
            $N_7=\begin{pmatrix}
            1999999999500&-826823851676&-826818852176&-582119098426&-835761802609&-582119098447&835761802637\\
-826823851676&2035761802334&1035761802584&-591057049402&591062048917&826823851673&-591057049411\\
-826818852176&1035761802584&2035761802334&-591062048902&591057049416&826818852173&-591062048911\\
-582119098426&-591057049402&-591062048902&1999999999500&835761802598&582119098437&582119098437\\
-835761802609&591062048917&591057049416&835761802598&1999999999500&835761802608&-582119098476\\
-582119098447&826823851673&826818852173&582119098437&835761802608&1999999999500&582119098416\\
835761802637&-591057049411&-591062048911&582119098437&-582119098476&582119098416&1999999999500\\
            \end{pmatrix}.$}
            \]

A routine computation shows that both $G_7$ and $G_7-\frac{1}{5}I_7$ are positive definite, where $I_n$ denotes the $n\times n$ identity matrix. Suppose that ${\bf z}\in\mathbb{Z}^7\setminus\{{\bf o}\}$ satisfies $Q_7({\bf z})<2$. Since
$${\bf z}^{\mathsf T}\Big(G_7-\frac{1}{5}I_7\Big){\bf z}>0,$$
we have
$$0<||{\bf z}||^2<5{\bf z}^{\mathsf T}G_7{\bf z}=5Q_7({\bf z})<10.$$
Enumeration by computer shows that
$$Z_7=\big\{{\bf z}\in\mathbb{Z}^7:\ 0<||{\bf z}||^2<10\big\}$$
contains $12434$ vectors. Computing $Q_7({\bf z})$ for all ${\bf z}\in Z_7$ by computer gives
$$\min_{{\bf z}\in Z_7}Q_7({\bf z})=Q_7\big((0,0,0,1,-1,0,-1)\big)=2,$$
which means
$$\min_{{\bf z}\in\mathbb{Z}^7\setminus\{{\bf o}\}}Q_7({\bf z})=2.$$
In other words, the shortest nonzero vectors of the corresponding
lattice $\Lambda$ have length $\sqrt2$ and
$$r(\Lambda)=1/\sqrt2.$$

The same inequality shows that every nonzero vector ${\bf z}$ satisfying $Q_{7}({\bf z})<4$ satisfies
$$||{\bf z}||^2<20.$$
We enumerate these vectors by computer and retain those satisfying $Q_7({\bf z})<4$. This yields exactly $127$ pairs of opposite vectors $\{\pm{\bf z}_1,\ldots,\pm{\bf z}_{127}\}$ such that
$$2\leq Q_7({\bf z}_1),\ldots,Q_7({\bf z}_{127})<4,$$
where
\begingroup
\small
\setlength{\LTleft}{\fill}
\setlength{\LTright}{\fill}
\begin{longtable}{@{}l@{\hspace{1.8em}}l@{\hspace{1.8em}}l@{}}
\(\mathbf{z}_{1}=(0,0,0,0,0,0,1)\), & \(\mathbf{z}_{2}=(0,0,0,0,0,1,-1)\), & \(\mathbf{z}_{3}=(0,0,0,0,0,1,0)\),\\
\(\mathbf{z}_{4}=(0,0,0,0,1,-1,0)\), & \(\mathbf{z}_{5}=(0,0,0,0,1,-1,1)\), & \(\mathbf{z}_{6}=(0,0,0,0,1,0,0)\),\\
\(\mathbf{z}_{7}=(0,0,0,0,1,0,1)\), & \(\mathbf{z}_{8}=(0,0,0,1,-2,1,-1)\), & \(\mathbf{z}_{9}=(0,0,0,1,-1,0,-1)\),\\
\(\mathbf{z}_{10}=(0,0,0,1,-1,0,0)\), & \(\mathbf{z}_{11}=(0,0,0,1,-1,1,-1)\), & \(\mathbf{z}_{12}=(0,0,0,1,-1,1,0)\),\\
\(\mathbf{z}_{13}=(0,0,0,1,0,-1,0)\), & \(\mathbf{z}_{14}=(0,0,0,1,0,0,-1)\), & \(\mathbf{z}_{15}=(0,0,0,1,0,0,0)\),\\
\(\mathbf{z}_{16}=(0,0,1,-1,1,-1,1)\), & \(\mathbf{z}_{17}=(0,0,1,0,-1,0,0)\), & \(\mathbf{z}_{18}=(0,0,1,0,0,-1,0)\),\\
\(\mathbf{z}_{19}=(0,0,1,0,0,-1,1)\), & \(\mathbf{z}_{20}=(0,0,1,0,0,0,0)\), & \(\mathbf{z}_{21}=(0,0,1,0,0,0,1)\),\\
\(\mathbf{z}_{22}=(0,0,1,0,1,-2,1)\), & \(\mathbf{z}_{23}=(0,0,1,0,1,-1,0)\), & \(\mathbf{z}_{24}=(0,0,1,0,1,-1,1)\),\\
\(\mathbf{z}_{25}=(0,0,1,1,-1,-1,0)\), & \(\mathbf{z}_{26}=(0,0,1,1,-1,0,-1)\), & \(\mathbf{z}_{27}=(0,0,1,1,-1,0,0)\),\\
\(\mathbf{z}_{28}=(0,0,1,1,0,-1,0)\), & \(\mathbf{z}_{29}=(0,0,1,1,0,-1,1)\), & \(\mathbf{z}_{30}=(0,0,1,1,0,0,0)\),\\
\(\mathbf{z}_{31}=(0,1,-1,-1,0,0,0)\), & \(\mathbf{z}_{32}=(0,1,-1,0,-1,0,0)\), & \(\mathbf{z}_{33}=(0,1,-1,0,-1,1,-1)\),\\
\(\mathbf{z}_{34}=(0,1,-1,0,0,-1,0)\), & \(\mathbf{z}_{35}=(0,1,-1,0,0,0,-1)\), & \(\mathbf{z}_{36}=(0,1,-1,0,0,0,0)\),\\
\(\mathbf{z}_{37}=(0,1,-1,1,-1,0,-1)\), & \(\mathbf{z}_{38}=(0,1,0,-1,1,-1,1)\), & \(\mathbf{z}_{39}=(0,1,0,0,-1,0,0)\),\\
\(\mathbf{z}_{40}=(0,1,0,0,0,-1,0)\), & \(\mathbf{z}_{41}=(0,1,0,0,0,-1,1)\), & \(\mathbf{z}_{42}=(0,1,0,0,0,0,0)\),\\
\(\mathbf{z}_{43}=(0,1,0,0,0,0,1)\), & \(\mathbf{z}_{44}=(0,1,0,0,1,-2,1)\), & \(\mathbf{z}_{45}=(0,1,0,0,1,-1,0)\),\\
\(\mathbf{z}_{46}=(0,1,0,0,1,-1,1)\), & \(\mathbf{z}_{47}=(0,1,0,1,-1,-1,0)\), & \(\mathbf{z}_{48}=(0,1,0,1,-1,0,-1)\),\\
\(\mathbf{z}_{49}=(0,1,0,1,-1,0,0)\), & \(\mathbf{z}_{50}=(0,1,0,1,0,-1,0)\), & \(\mathbf{z}_{51}=(0,1,0,1,0,-1,1)\),\\
\(\mathbf{z}_{52}=(0,1,0,1,0,0,0)\), & \(\mathbf{z}_{53}=(0,1,1,0,0,-1,1)\), & \(\mathbf{z}_{54}=(0,1,1,1,-1,-1,0)\),\\
\(\mathbf{z}_{55}=(0,1,1,1,-1,-1,1)\), & \(\mathbf{z}_{56}=(0,1,1,1,-1,0,0)\), & \(\mathbf{z}_{57}=(0,1,1,1,0,-2,1)\),\\
\(\mathbf{z}_{58}=(0,1,1,1,0,-1,0)\), & \(\mathbf{z}_{59}=(0,1,1,1,0,-1,1)\), & \(\mathbf{z}_{60}=(0,1,1,2,-1,-1,0)\),\\
\(\mathbf{z}_{61}=(1,-1,0,0,0,1,-1)\), & \(\mathbf{z}_{62}=(1,0,-1,0,0,1,-1)\), & \(\mathbf{z}_{63}=(1,0,0,-1,1,0,0)\),\\
\(\mathbf{z}_{64}=(1,0,0,0,-1,1,-1)\), & \(\mathbf{z}_{65}=(1,0,0,0,0,0,-1)\), & \(\mathbf{z}_{66}=(1,0,0,0,0,0,0)\),\\
\(\mathbf{z}_{67}=(1,0,0,0,0,1,-1)\), & \(\mathbf{z}_{68}=(1,0,0,0,0,1,0)\), & \(\mathbf{z}_{69}=(1,0,0,0,1,-1,0)\),\\
\(\mathbf{z}_{70}=(1,0,0,0,1,0,-1)\), & \(\mathbf{z}_{71}=(1,0,0,0,1,0,0)\), & \(\mathbf{z}_{72}=(1,0,0,1,-1,0,-1)\),\\
\(\mathbf{z}_{73}=(1,0,0,1,-1,1,-2)\), & \(\mathbf{z}_{74}=(1,0,0,1,-1,1,-1)\), & \(\mathbf{z}_{75}=(1,0,0,1,0,0,-1)\),\\
\(\mathbf{z}_{76}=(1,0,0,1,0,0,0)\), & \(\mathbf{z}_{77}=(1,0,0,1,0,1,-1)\), & \(\mathbf{z}_{78}=(1,0,1,0,0,-1,0)\),\\
\(\mathbf{z}_{79}=(1,0,1,0,0,0,-1)\), & \(\mathbf{z}_{80}=(1,0,1,0,0,0,0)\), & \(\mathbf{z}_{81}=(1,0,1,0,1,-1,0)\),\\
\(\mathbf{z}_{82}=(1,0,1,0,1,-1,1)\), & \(\mathbf{z}_{83}=(1,0,1,0,1,0,0)\), & \(\mathbf{z}_{84}=(1,0,1,1,-1,0,-1)\),\\
\(\mathbf{z}_{85}=(1,0,1,1,-1,0,0)\), & \(\mathbf{z}_{86}=(1,0,1,1,-1,1,-1)\), & \(\mathbf{z}_{87}=(1,0,1,1,0,-1,-1)\),\\
\(\mathbf{z}_{88}=(1,0,1,1,0,-1,0)\), & \(\mathbf{z}_{89}=(1,0,1,1,0,0,-1)\), & \(\mathbf{z}_{90}=(1,0,1,1,0,0,0)\),\\
\(\mathbf{z}_{91}=(1,0,1,1,1,-1,0)\), & \(\mathbf{z}_{92}=(1,0,1,2,-1,0,-1)\), & \(\mathbf{z}_{93}=(1,1,-1,0,0,0,0)\),\\
\(\mathbf{z}_{94}=(1,1,-1,0,0,1,-1)\), & \(\mathbf{z}_{95}=(1,1,-1,1,-1,1,-1)\), & \(\mathbf{z}_{96}=(1,1,-1,1,0,0,-1)\),\\
\(\mathbf{z}_{97}=(1,1,0,0,0,-1,0)\), & \(\mathbf{z}_{98}=(1,1,0,0,0,0,-1)\), & \(\mathbf{z}_{99}=(1,1,0,0,0,0,0)\),\\
\(\mathbf{z}_{100}=(1,1,0,0,1,-1,0)\), & \(\mathbf{z}_{101}=(1,1,0,0,1,-1,1)\), & \(\mathbf{z}_{102}=(1,1,0,0,1,0,0)\),\\
\(\mathbf{z}_{103}=(1,1,0,1,-1,0,-1)\), & \(\mathbf{z}_{104}=(1,1,0,1,-1,0,0)\), & \(\mathbf{z}_{105}=(1,1,0,1,-1,1,-1)\),\\
\(\mathbf{z}_{106}=(1,1,0,1,0,-1,-1)\), & \(\mathbf{z}_{107}=(1,1,0,1,0,-1,0)\), & \(\mathbf{z}_{108}=(1,1,0,1,0,0,-1)\),\\
\(\mathbf{z}_{109}=(1,1,0,1,0,0,0)\), & \(\mathbf{z}_{110}=(1,1,0,1,1,-1,0)\), & \(\mathbf{z}_{111}=(1,1,0,2,-1,0,-1)\),\\
\(\mathbf{z}_{112}=(1,1,1,1,-1,-1,0)\), & \(\mathbf{z}_{113}=(1,1,1,1,-1,0,-1)\), & \(\mathbf{z}_{114}=(1,1,1,1,-1,0,0)\),\\
\(\mathbf{z}_{115}=(1,1,1,1,0,-1,0)\), & \(\mathbf{z}_{116}=(1,1,1,1,0,-1,1)\), & \(\mathbf{z}_{117}=(1,1,1,1,0,0,0)\),\\
\(\mathbf{z}_{118}=(1,1,1,2,-2,0,-1)\), & \(\mathbf{z}_{119}=(1,1,1,2,-1,-1,-1)\), & \(\mathbf{z}_{120}=(1,1,1,2,-1,-1,0)\),\\
\(\mathbf{z}_{121}=(1,1,1,2,-1,0,-1)\), & \(\mathbf{z}_{122}=(1,1,1,2,-1,0,0)\), & \(\mathbf{z}_{123}=(1,1,1,2,0,-1,0)\),\\
\(\mathbf{z}_{124}=(1,2,0,1,0,-1,0)\), & \(\mathbf{z}_{125}=(2,0,1,1,0,0,-1)\), & \(\mathbf{z}_{126}=(2,1,0,1,0,0,-1)\),\\
\(\mathbf{z}_{127}=(2,1,1,2,-1,0,-1)\). & & \\
\end{longtable}
\endgroup
Therefore, we have
$$\big\{{\bf z}\in \mathbb{Z}^7:\ 2\leq Q_7({\bf z})<4\big\}=\{\pm{\bf z}_1,\ldots,\pm{\bf z}_{127}\}.$$

Using Remark 4 and the algorithm described in \cite{Av01}, we computed all vertices of the Voronoi cell $V(\Lambda)$ of the lattice $\Lambda$ corresponding to $Q_7({\bf z})$ by computer. The computation shows that $V(\Lambda)$ has $254$ facets and $40320$ vertices, and that the maximum distance from the origin ${\bf o}$ to a vertex is
\[
\resizebox{1.002\linewidth}{!}{%
$\displaystyle
\begin{aligned}
 R(\Lambda)&=\left(\frac{400649880198127881308780272663122998605428023002937746204041887507637863065104887397836871724}{348734693306100951459852010372216453497786214276229048810574846905601413447280352881430974187}\right)^{\frac{1}{2}}\\[0.6em]
&=1.0718522681573623\cdots,
\end{aligned}
$%
}
\]
which shows that
$$\gamma^*(B^7)\leq R(\Lambda)/r(\Lambda)=\sqrt2R(\Lambda)=1.51582801\cdots<\sqrt{7/3}.$$

Enumeration by computer shows that all the $Q_7({\bf z})$ satisfying $2\leq Q_7({\bf z})<4$ fall into the five intervals listed in Table 3.

\begin{table}[ht]
\newcommand{\AD}[1]{{\begin{tabular}{@{}c@{}}#1\end{tabular}}}
\renewcommand{\arraystretch}{1.2}
\begin{tabular}{|c|c|c|c|c|c|}
\hline ${\rm Interval}\ I$ & $[2.0,2.1)$ & $[2.3,2.4)$ & $[2.8,3.0)$ & $[3.3,3.4)$ & $[3.7,4.0)$\\
\hline ${\rm card}\big\{{\bf z}\in\mathbb{Z}^7:\ Q_{7}({\bf z})\in I\big\}$ & $48$ & $30$ & $70$ & $20$ & $86$\\
\hline
\end{tabular}

\vspace{0.3cm}
\centerline{Table 3. The distribution of small values of $Q_7({\bf z})$.}
\end{table}

As a comparison, we recall (see \cite[p.125]{Co01}) that the first three nonzero shells of the positive definite quadratic form $Q_{{\scriptscriptstyle\sqrt{4/3}}E^*_{7}}({\bf z})$ corresponding to $\sqrt{4/3}E^*_{7}$ are listed in Table 4 as following.

\begin{table}[ht]
\newcommand{\AD}[1]{{\begin{tabular}{@{}c@{}}#1\end{tabular}}}
\renewcommand{\arraystretch}{1.3}
\begin{tabular}{|c|c|c|c|}
\hline ${\rm Norm\ } N$ & $2$ & $8/3$ & $14/3$\\
\hline ${\rm card}\big\{{\bf z}\in\mathbb{Z}^7:\ Q_{{\scriptscriptstyle\sqrt{4/3}}E^*_{7}}({\bf z})=N\big\}$ & $56$ & $126$ & $576$\\
\hline
\end{tabular}

\vspace{0.3cm}
\centerline{Table 4. The distribution of small values of $Q_{{\scriptscriptstyle\sqrt{4/3}}E^*_{7}}({\bf z})$.}
\end{table}

\bigskip
\noindent
{\Large\bf 3.2. The Eight-Dimensional Case}

\medskip
Let $Q_8({\bf z})={\bf z}^{\mathsf T}G_8{\bf z}$, where $G_8=N_8/d_8$, $d_8=7999999995$ and
\[
            \resizebox{\textwidth}{!}{%
            $N_8=\begin{pmatrix}
            16000000008&-7991363624&21438716&-21408352&-23669032&-18641312&-11142664&14946528\\
-7991363624&16000000000&-7999845964&-4281740&-8469568&30643280&-30526408&8218328\\
21438716&-7999845964&16000000000&-7995872292&14921716&-3411288&-18657172&-8019003720\\
-21408352&-4281740&-7995872292&15999999998&-8010744664&10744660&2982504&6540884\\
-23669032&-8469568&14921716&-8010744664&16038589600&-8019294800&23305344&-29705736\\
-18641312&30643280&-3411288&10744660&-8019294800&16000000000&-7966031272&7526720\\
-11142664&-30526408&-18657172&2982504&23305344&-7966031272&16000058256&814576\\
14946528&8218328&-8019003720&6540884&-29705736&7526720&814576&16038007440\\
            \end{pmatrix}.$}
            \]

By the definition of the lattice $E_8$ (see \cite[p.120]{Co01}), a Gram matrix $G_{E_8}$ of the lattice $E_8$ is
\[
G_{E_{8}}=\begin{pmatrix}
2&-1&0&0&0&0&0&0\\
-1&2&-1&0&0&0&0&0\\
0&-1&2&-1&0&0&0&-1\\
0&0&-1&2&-1&0&0&0\\
0&0&0&-1&2&-1&0&0\\
0&0&0&0&-1&2&-1&0\\
0&0&0&0&0&-1&2&0\\
0&0&-1&0&0&0&0&2\\
\end{pmatrix}.
\]
A routine computation shows that both $G_8$ and $G_8-\frac{2}{3}G_{E_8}$ are positive definite. Let $Q_{E_8}({\bf z})$ denote the positive definite quadratic form corresponding to $E_8$. Suppose that ${\bf z}\in\mathbb{Z}^8\setminus\{{\bf o}\}$ satisfies $Q_8({\bf z})<2$. Since
$${\bf z}^{\mathsf T}\Big(G_8-\frac{2}{3}G_{E_8}\Big){\bf z}>0,$$
we have
$$0<Q_{E_8}({\bf z})={\bf z}^{\mathsf T}G_{E_8}{\bf z}<\frac32{\bf z}^{\mathsf T}G_8{\bf z}=\frac32 Q_8({\bf z})<3.$$
It follows from the theta series of $E_8$ given in \cite[p.123]{Co01} that the set
$$Z_8=\big\{{\bf z}\in\mathbb Z^8:\ 0<Q_{E_8}({\bf z})<3\big\}
     =\big\{{\bf z}\in\mathbb Z^8:\ Q_{E_8}({\bf z})=2\big\}$$
contains \(240\) vectors. Computing $Q_8({\bf z})$ for all ${\bf z}\in Z_8$ by computer gives
$$\min_{{\bf z}\in Z_8}Q_8({\bf z})=Q_8\bigl((0,0,0,1,1,1,0,0)\bigr)=2,$$
which means
$$\min_{{\bf z}\in\mathbb Z^8\setminus\{{\bf o}\}}Q_8(\mathbf z)=2.$$
In other words, the shortest nonzero vectors of the corresponding
lattice $\Lambda$ have length $\sqrt2$ and
$$r(\Lambda)=1/\sqrt2.$$

The same inequality shows that every nonzero vector ${\bf z}$ satisfying $Q_{8}({\bf z})<4$ satisfies
$$Q_{E_8}({\bf z})<6.$$
By the theta series of $E_8$ given in \cite[p.123]{Co01}, this condition implies
$$Q_{E_8}({\bf z})=2\  {\rm or}\  4.$$
We enumerate these vectors by computer and retain those satisfying $Q_8({\bf z})<4$. This yields exactly $255$ pairs of opposite vectors $\{\pm\mathbf z_1,\ldots,\pm\mathbf z_{255}\}$ such that
$$2\le Q_8({\bf z}_1),\ldots,Q_8({\bf z}_{255})<4,$$
where
\begingroup
\small
\setlength{\LTleft}{\fill}
\setlength{\LTright}{\fill}
\begin{longtable}{@{}l@{\hspace{1.2em}}l@{\hspace{1.2em}}l@{}}
\(\mathbf{z}_{1}=(1,1,2,2,1,1,1,1)\), & \(\mathbf{z}_{2}=(2,3,4,4,3,2,2,2)\), & \(\mathbf{z}_{3}=(1,2,2,2,2,1,1,1)\), \\
\(\mathbf{z}_{4}=(1,2,3,2,2,1,1,2)\), & \(\mathbf{z}_{5}=(1,1,1,2,1,1,1,0)\), & \(\mathbf{z}_{6}=(0,0,1,0,0,0,0,1)\), \\
\(\mathbf{z}_{7}=(0,0,0,0,0,0,0,1)\), & \(\mathbf{z}_{8}=(1,2,3,2,2,1,1,1)\), & \(\mathbf{z}_{9}=(1,1,1,2,1,1,1,1)\), \\
\(\mathbf{z}_{10}=(0,0,1,0,0,0,0,0)\), & \(\mathbf{z}_{11}=(0,0,2,1,0,0,0,1)\), & \(\mathbf{z}_{12}=(1,1,0,1,1,1,1,0)\), \\
\(\mathbf{z}_{13}=(1,2,4,3,2,1,1,2)\), & \(\mathbf{z}_{14}=(0,1,2,1,1,0,0,1)\), & \(\mathbf{z}_{15}=(0,0,1,1,0,0,0,0)\), \\
\(\mathbf{z}_{16}=(1,1,1,1,1,1,1,1)\), & \(\mathbf{z}_{17}=(1,2,3,3,2,1,1,1)\), & \(\mathbf{z}_{18}=(0,1,1,1,1,0,0,0)\), \\
\(\mathbf{z}_{19}=(0,0,1,1,0,0,0,1)\), & \(\mathbf{z}_{20}=(1,1,1,1,1,1,1,0)\), & \(\mathbf{z}_{21}=(1,2,3,3,2,1,1,2)\), \\
\(\mathbf{z}_{22}=(0,1,1,1,1,0,0,1)\), & \(\mathbf{z}_{23}=(0,0,0,1,0,0,0,0)\), & \(\mathbf{z}_{24}=(1,1,2,1,1,1,1,1)\), \\
\(\mathbf{z}_{25}=(1,2,2,3,2,1,1,1)\), & \(\mathbf{z}_{26}=(0,1,0,1,1,0,0,0)\), & \(\mathbf{z}_{27}=(0,1,3,2,1,0,0,2)\), \\
\(\mathbf{z}_{28}=(1,2,4,4,2,1,1,2)\), & \(\mathbf{z}_{29}=(0,1,2,2,1,0,0,1)\), & \(\mathbf{z}_{30}=(1,3,4,4,3,1,1,2)\), \\
\(\mathbf{z}_{31}=(0,1,1,2,1,0,0,1)\), & \(\mathbf{z}_{32}=(1,3,4,4,3,2,1,2)\), & \(\mathbf{z}_{33}=(1,0,0,0,0,0,1,0)\), \\
\(\mathbf{z}_{34}=(0,1,2,2,1,1,0,1)\), & \(\mathbf{z}_{35}=(1,3,5,4,3,2,1,3)\), & \(\mathbf{z}_{36}=(1,3,5,4,3,2,1,2)\), \\
\(\mathbf{z}_{37}=(1,1,2,1,1,0,1,1)\), & \(\mathbf{z}_{38}=(0,1,1,1,1,1,0,1)\), & \(\mathbf{z}_{39}=(1,2,3,3,2,2,1,2)\), \\
\(\mathbf{z}_{40}=(1,1,1,1,1,0,1,0)\), & \(\mathbf{z}_{41}=(0,1,1,1,1,1,0,0)\), & \(\mathbf{z}_{42}=(1,2,3,3,2,2,1,1)\), \\
\(\mathbf{z}_{43}=(1,1,1,1,1,0,1,1)\), & \(\mathbf{z}_{44}=(1,3,4,3,3,2,1,2)\), & \(\mathbf{z}_{45}=(2,4,6,5,4,3,2,3)\), \\
\(\mathbf{z}_{46}=(0,1,2,1,1,1,0,1)\), & \(\mathbf{z}_{47}=(1,2,4,3,2,2,1,2)\), & \(\mathbf{z}_{48}=(0,0,1,0,0,-1,0,1)\), \\
\(\mathbf{z}_{49}=(1,2,2,2,2,2,1,1)\), & \(\mathbf{z}_{50}=(1,1,2,2,1,0,1,1)\), & \(\mathbf{z}_{51}=(0,0,0,0,0,1,0,0)\), \\
\(\mathbf{z}_{52}=(1,2,3,2,2,2,1,2)\), & \(\mathbf{z}_{53}=(1,2,3,2,2,2,1,1)\), & \(\mathbf{z}_{54}=(0,0,1,0,0,1,0,0)\), \\
\(\mathbf{z}_{55}=(0,0,1,1,0,-1,0,0)\), & \(\mathbf{z}_{56}=(0,0,1,1,0,-1,0,1)\), & \(\mathbf{z}_{57}=(0,0,0,1,0,-1,0,0)\), \\
\(\mathbf{z}_{58}=(1,1,2,2,2,2,0,1)\), & \(\mathbf{z}_{59}=(1,0,0,0,1,1,0,0)\), & \(\mathbf{z}_{60}=(1,0,1,1,1,1,0,0)\), \\
\(\mathbf{z}_{61}=(0,1,1,1,0,0,1,1)\), & \(\mathbf{z}_{62}=(1,0,1,1,1,1,0,1)\), & \(\mathbf{z}_{63}=(0,1,1,1,0,0,1,0)\), \\
\(\mathbf{z}_{64}=(2,2,3,3,3,2,1,2)\), & \(\mathbf{z}_{65}=(1,1,1,1,2,1,0,1)\), & \(\mathbf{z}_{66}=(1,0,0,1,1,1,0,0)\), \\
\(\mathbf{z}_{67}=(0,1,2,1,0,0,1,1)\), & \(\mathbf{z}_{68}=(0,1,0,0,0,0,1,0)\), & \(\mathbf{z}_{69}=(1,1,2,2,2,1,0,1)\), \\
\(\mathbf{z}_{70}=(0,0,0,0,1,0,-1,0)\), & \(\mathbf{z}_{71}=(0,1,1,0,0,0,1,1)\), & \(\mathbf{z}_{72}=(1,1,1,2,2,1,0,0)\), \\
\(\mathbf{z}_{73}=(0,1,1,0,0,0,1,0)\), & \(\mathbf{z}_{74}=(1,1,1,2,2,1,0,1)\), & \(\mathbf{z}_{75}=(0,0,1,1,1,0,-1,1)\), \\
\(\mathbf{z}_{76}=(1,1,2,3,2,1,0,1)\), & \(\mathbf{z}_{77}=(0,0,0,1,1,0,-1,0)\), & \(\mathbf{z}_{78}=(1,0,0,0,1,0,0,0)\), \\
\(\mathbf{z}_{79}=(1,0,1,1,1,0,0,1)\), & \(\mathbf{z}_{80}=(1,0,0,1,1,0,0,0)\), & \(\mathbf{z}_{81}=(1,3,4,3,2,2,1,2)\), \\
\(\mathbf{z}_{82}=(0,0,0,1,1,0,0,0)\), & \(\mathbf{z}_{83}=(1,1,2,3,2,1,1,1)\), & \(\mathbf{z}_{84}=(1,2,2,1,1,1,1,1)\), \\
\(\mathbf{z}_{85}=(0,0,1,1,1,0,0,1)\), & \(\mathbf{z}_{86}=(1,2,1,1,1,1,1,0)\), & \(\mathbf{z}_{87}=(0,0,1,1,1,0,0,0)\), \\
\(\mathbf{z}_{88}=(0,1,2,0,0,0,0,1)\), & \(\mathbf{z}_{89}=(1,1,1,2,2,1,1,1)\), & \(\mathbf{z}_{90}=(1,2,3,2,1,1,1,1)\), \\
\(\mathbf{z}_{91}=(0,1,1,0,0,0,0,0)\), & \(\mathbf{z}_{92}=(1,2,2,2,1,1,1,0)\), & \(\mathbf{z}_{93}=(1,2,3,2,1,1,1,2)\), \\
\(\mathbf{z}_{94}=(0,1,1,0,0,0,0,1)\), & \(\mathbf{z}_{95}=(0,0,0,0,1,0,0,0)\), & \(\mathbf{z}_{96}=(1,1,2,2,2,1,1,1)\), \\
\(\mathbf{z}_{97}=(1,2,2,2,1,1,1,1)\), & \(\mathbf{z}_{98}=(0,1,0,0,0,0,0,0)\), & \(\mathbf{z}_{99}=(0,1,2,1,0,0,0,1)\), \\
\(\mathbf{z}_{100}=(1,0,0,1,1,1,1,0)\), & \(\mathbf{z}_{101}=(1,3,4,3,2,1,1,2)\), & \(\mathbf{z}_{102}=(0,2,2,1,1,0,0,1)\), \\
\(\mathbf{z}_{103}=(1,2,3,3,1,1,1,1)\), & \(\mathbf{z}_{104}=(0,1,1,1,0,0,0,0)\), & \(\mathbf{z}_{105}=(1,3,3,3,2,1,1,1)\), \\
\(\mathbf{z}_{106}=(0,1,1,1,0,0,0,1)\), & \(\mathbf{z}_{107}=(1,3,3,3,2,1,1,2)\), & \(\mathbf{z}_{108}=(0,1,0,1,0,0,0,0)\), \\
\(\mathbf{z}_{109}=(1,3,5,4,2,1,1,2)\), & \(\mathbf{z}_{110}=(0,2,3,2,1,0,0,1)\), & \(\mathbf{z}_{111}=(0,2,3,2,1,0,0,2)\), \\
\(\mathbf{z}_{112}=(1,3,4,4,2,1,1,2)\), & \(\mathbf{z}_{113}=(0,2,2,2,1,0,0,1)\), & \(\mathbf{z}_{114}=(1,1,1,0,0,0,1,0)\), \\
\(\mathbf{z}_{115}=(1,1,1,0,0,0,1,1)\), & \(\mathbf{z}_{116}=(0,1,2,2,2,1,0,1)\), & \(\mathbf{z}_{117}=(1,2,4,4,3,2,1,2)\), \\
\(\mathbf{z}_{118}=(1,1,0,0,0,0,1,0)\), & \(\mathbf{z}_{119}=(1,1,2,1,0,0,1,1)\), & \(\mathbf{z}_{120}=(0,0,0,1,1,1,0,0)\), \\
\(\mathbf{z}_{121}=(1,2,2,1,1,0,1,1)\), & \(\mathbf{z}_{122}=(1,2,3,3,3,2,1,2)\), & \(\mathbf{z}_{123}=(1,1,1,1,0,0,1,0)\), \\
\(\mathbf{z}_{124}=(0,0,1,1,1,1,0,1)\), & \(\mathbf{z}_{125}=(1,2,3,3,3,2,1,1)\), & \(\mathbf{z}_{126}=(1,1,1,1,0,0,1,1)\), \\
\(\mathbf{z}_{127}=(0,0,1,1,1,1,0,0)\), & \(\mathbf{z}_{128}=(1,2,4,3,3,2,1,2)\), & \(\mathbf{z}_{129}=(0,0,2,1,1,1,0,1)\), \\
\(\mathbf{z}_{130}=(0,0,-1,0,1,1,0,0)\), & \(\mathbf{z}_{131}=(1,2,3,2,1,0,1,1)\), & \(\mathbf{z}_{132}=(0,1,1,0,0,-1,0,0)\), \\
\(\mathbf{z}_{133}=(0,0,0,0,1,1,0,1)\), & \(\mathbf{z}_{134}=(1,2,3,2,1,0,1,2)\), & \(\mathbf{z}_{135}=(0,1,1,0,0,-1,0,1)\), \\
\(\mathbf{z}_{136}=(0,0,0,0,1,1,0,0)\), & \(\mathbf{z}_{137}=(1,1,2,2,2,2,1,1)\), & \(\mathbf{z}_{138}=(1,2,2,2,1,0,1,1)\), \\
\(\mathbf{z}_{139}=(0,1,0,0,0,-1,0,0)\), & \(\mathbf{z}_{140}=(0,0,1,0,1,1,0,1)\), & \(\mathbf{z}_{141}=(0,1,2,1,0,-1,0,1)\), \\
\(\mathbf{z}_{142}=(0,1,1,1,0,-1,0,0)\), & \(\mathbf{z}_{143}=(0,1,1,1,0,-1,0,1)\), & \(\mathbf{z}_{144}=(2,4,6,5,4,3,1,3)\), \\
\(\mathbf{z}_{145}=(1,0,0,0,0,1,0,0)\), & \(\mathbf{z}_{146}=(0,1,2,2,1,0,1,1)\), & \(\mathbf{z}_{147}=(2,2,2,2,2,2,1,1)\), \\
\(\mathbf{z}_{148}=(1,1,0,0,1,1,0,0)\), & \(\mathbf{z}_{149}=(1,1,2,1,1,1,0,1)\), & \(\mathbf{z}_{150}=(2,3,4,3,3,2,1,2)\), \\
\(\mathbf{z}_{151}=(0,1,1,1,1,0,1,1)\), & \(\mathbf{z}_{152}=(1,1,1,1,1,1,0,0)\), & \(\mathbf{z}_{153}=(0,0,1,1,0,0,1,1)\), \\
\(\mathbf{z}_{154}=(0,1,1,1,1,0,1,0)\), & \(\mathbf{z}_{155}=(1,1,1,1,1,1,0,1)\), & \(\mathbf{z}_{156}=(0,0,1,1,0,0,1,0)\), \\
\(\mathbf{z}_{157}=(0,1,2,1,1,0,1,1)\), & \(\mathbf{z}_{158}=(1,1,0,1,1,1,0,0)\), & \(\mathbf{z}_{159}=(2,3,5,4,3,2,1,2)\), \\
\(\mathbf{z}_{160}=(1,2,3,2,2,1,0,1)\), & \(\mathbf{z}_{161}=(2,3,5,4,3,2,1,3)\), & \(\mathbf{z}_{162}=(1,2,3,2,2,1,0,2)\), \\
\(\mathbf{z}_{163}=(1,1,2,2,1,1,0,1)\), & \(\mathbf{z}_{164}=(0,0,0,0,0,0,1,0)\), & \(\mathbf{z}_{165}=(2,3,4,4,3,2,1,2)\), \\
\(\mathbf{z}_{166}=(1,2,2,2,2,1,0,1)\), & \(\mathbf{z}_{167}=(1,1,1,2,1,1,0,0)\), & \(\mathbf{z}_{168}=(0,0,1,0,0,0,1,1)\), \\
\(\mathbf{z}_{169}=(1,1,1,2,1,1,0,1)\), & \(\mathbf{z}_{170}=(0,0,1,0,0,0,1,0)\), & \(\mathbf{z}_{171}=(1,2,4,3,2,1,0,2)\), \\
\(\mathbf{z}_{172}=(2,4,6,5,4,2,1,3)\), & \(\mathbf{z}_{173}=(1,2,3,3,2,1,0,1)\), & \(\mathbf{z}_{174}=(1,2,3,3,2,1,0,2)\), \\
\(\mathbf{z}_{175}=(0,0,0,1,0,0,-1,0)\), & \(\mathbf{z}_{176}=(1,2,2,3,2,1,0,1)\), & \(\mathbf{z}_{177}=(-1,0,1,1,0,0,0,1)\), \\
\(\mathbf{z}_{178}=(2,2,3,2,2,1,1,1)\), & \(\mathbf{z}_{179}=(1,0,0,0,0,0,0,0)\), & \(\mathbf{z}_{180}=(0,1,2,2,1,1,1,1)\), \\
\(\mathbf{z}_{181}=(2,2,2,2,2,1,1,1)\), & \(\mathbf{z}_{182}=(1,1,0,0,1,0,0,0)\), & \(\mathbf{z}_{183}=(-1,0,1,0,0,0,0,1)\), \\
\(\mathbf{z}_{184}=(1,0,-1,0,0,0,0,0)\), & \(\mathbf{z}_{185}=(2,2,4,3,2,1,1,2)\), & \(\mathbf{z}_{186}=(1,1,2,1,1,0,0,1)\), \\
\(\mathbf{z}_{187}=(2,3,4,3,3,1,1,2)\), & \(\mathbf{z}_{188}=(0,1,1,1,1,1,1,1)\), & \(\mathbf{z}_{189}=(2,2,3,3,2,1,1,1)\), \\
\(\mathbf{z}_{190}=(1,1,1,1,1,0,0,0)\), & \(\mathbf{z}_{191}=(0,1,1,1,1,1,1,0)\), & \(\mathbf{z}_{192}=(2,2,3,3,2,1,1,2)\), \\
\(\mathbf{z}_{193}=(1,1,1,1,1,0,0,1)\), & \(\mathbf{z}_{194}=(1,0,0,1,0,0,0,0)\), & \(\mathbf{z}_{195}=(0,1,2,1,1,1,1,1)\), \\
\(\mathbf{z}_{196}=(1,1,0,1,1,0,0,0)\), & \(\mathbf{z}_{197}=(2,3,5,4,3,1,1,2)\), & \(\mathbf{z}_{198}=(2,3,5,4,3,1,1,3)\), \\
\(\mathbf{z}_{199}=(1,1,2,2,1,0,0,1)\), & \(\mathbf{z}_{200}=(0,0,0,0,0,1,1,0)\), & \(\mathbf{z}_{201}=(2,3,4,4,3,1,1,2)\), \\
\(\mathbf{z}_{202}=(1,2,2,2,2,0,0,1)\), & \(\mathbf{z}_{203}=(1,1,1,2,1,0,0,0)\), & \(\mathbf{z}_{204}=(1,1,1,2,1,0,0,1)\), \\
\(\mathbf{z}_{205}=(1,2,3,3,2,0,0,2)\), & \(\mathbf{z}_{206}=(1,2,2,1,0,0,1,1)\), & \(\mathbf{z}_{207}=(0,-1,0,1,1,1,0,0)\), \\
\(\mathbf{z}_{208}=(0,1,1,0,-1,-1,0,0)\), & \(\mathbf{z}_{209}=(0,1,1,0,-1,-1,0,1)\), & \(\mathbf{z}_{210}=(0,-1,0,0,1,1,0,0)\), \\
\(\mathbf{z}_{211}=(2,3,4,3,2,2,1,2)\), & \(\mathbf{z}_{212}=(1,2,2,1,1,1,0,1)\), & \(\mathbf{z}_{213}=(2,4,6,4,3,2,1,3)\), \\
\(\mathbf{z}_{214}=(1,2,3,2,1,1,0,1)\), & \(\mathbf{z}_{215}=(2,4,5,4,3,2,1,2)\), & \(\mathbf{z}_{216}=(1,2,3,2,1,1,0,2)\), \\
\(\mathbf{z}_{217}=(2,4,5,4,3,2,1,3)\), & \(\mathbf{z}_{218}=(1,2,2,2,1,1,0,1)\), & \(\mathbf{z}_{219}=(2,4,6,5,3,2,1,3)\), \\
\(\mathbf{z}_{220}=(1,3,4,3,2,1,0,2)\), & \(\mathbf{z}_{221}=(1,1,1,0,0,0,0,0)\), & \(\mathbf{z}_{222}=(2,3,3,2,2,1,1,1)\), \\
\(\mathbf{z}_{223}=(1,1,1,0,0,0,0,1)\), & \(\mathbf{z}_{224}=(0,1,2,2,2,1,1,1)\), & \(\mathbf{z}_{225}=(2,2,2,2,1,1,1,1)\), \\
\(\mathbf{z}_{226}=(1,1,0,0,0,0,0,0)\), & \(\mathbf{z}_{227}=(2,3,4,3,2,1,1,1)\), & \(\mathbf{z}_{228}=(1,1,2,1,0,0,0,1)\), \\
\(\mathbf{z}_{229}=(0,0,0,1,1,1,1,0)\), & \(\mathbf{z}_{230}=(2,3,4,3,2,1,1,2)\), & \(\mathbf{z}_{231}=(1,2,2,1,1,0,0,1)\), \\
\(\mathbf{z}_{232}=(1,1,1,1,0,0,0,0)\), & \(\mathbf{z}_{233}=(0,0,1,1,1,1,1,1)\), & \(\mathbf{z}_{234}=(2,3,3,3,2,1,1,1)\), \\
\(\mathbf{z}_{235}=(1,2,1,1,1,0,0,0)\), & \(\mathbf{z}_{236}=(1,1,1,1,0,0,0,1)\), & \(\mathbf{z}_{237}=(0,0,1,1,1,1,1,0)\), \\
\(\mathbf{z}_{238}=(2,3,3,3,2,1,1,2)\), & \(\mathbf{z}_{239}=(1,1,0,1,0,0,0,0)\), & \(\mathbf{z}_{240}=(2,3,5,4,2,1,1,2)\), \\
\(\mathbf{z}_{241}=(1,2,3,2,1,0,0,1)\), & \(\mathbf{z}_{242}=(2,4,5,4,3,1,1,2)\), & \(\mathbf{z}_{243}=(1,2,3,2,1,0,0,2)\), \\
\(\mathbf{z}_{244}=(2,4,5,4,3,1,1,3)\), & \(\mathbf{z}_{245}=(0,0,0,0,1,1,1,0)\), & \(\mathbf{z}_{246}=(2,3,4,4,2,1,1,2)\), \\
\(\mathbf{z}_{247}=(1,2,2,2,1,0,0,1)\), & \(\mathbf{z}_{248}=(2,4,6,5,3,1,1,3)\), & \(\mathbf{z}_{249}=(0,1,2,2,2,2,1,1)\), \\
\(\mathbf{z}_{250}=(1,1,0,0,0,-1,0,0)\), & \(\mathbf{z}_{251}=(1,1,1,1,0,-1,0,0)\), & \(\mathbf{z}_{252}=(2,2,2,2,2,1,0,1)\), \\
\(\mathbf{z}_{253}=(1,2,2,1,0,0,0,1)\), & \(\mathbf{z}_{254}=(2,4,5,4,2,1,1,2)\), & \(\mathbf{z}_{255}=(2,3,4,3,2,1,0,2)\). \\
\end{longtable}
\endgroup

Therefore, we have
$$\big\{{\bf z}\in \mathbb{Z}^8:\ 2\leq Q_8({\bf z})<4\big\}=\{\pm{\bf z}_1,\ldots,\pm{\bf z}_{255}\}.$$

Using Remark 4 and the algorithm described in \cite{Av01}, we computed all vertices of the Voronoi cell $V(\Lambda)$ of the lattice $\Lambda$ corresponding to $Q_8({\bf z})$ by computer. The computation shows that $V(\Lambda)$ has $510$ facets and $172800$ vertices, and that the maximum distance from the origin ${\bf o}$ to a vertex is
\[
\resizebox{0.99\linewidth}{!}{%
$\displaystyle
\begin{aligned}
 R(\Lambda)&=\left(\frac{265716457506397653493171307697740018550559863089465523139020783198810022846440906787}{265747715525373181682796646403850821917661505831567269814397604656016490008569878000}\right)^{\frac{1}{2}}\\[0.6em]
&=0.9999411868167285\cdots,
\end{aligned}
$%
}
\]
which shows that
$$\gamma^*(B^8)\leq R(\Lambda)/r(\Lambda)=\sqrt2R(\Lambda)=1.41413038\cdots<\sqrt2.$$

Enumeration by computer shows that all the $Q_8({\bf z})$ satisfying $2\leq Q_8({\bf z})<4.1$ fall into the three intervals listed in Table 5.

\begin{table}[ht]
\vspace{-0.06cm}
\newcommand{\AD}[1]{{\begin{tabular}{@{}c@{}}#1\end{tabular}}}
\renewcommand{\arraystretch}{1.2}
\begin{tabular}{|c|c|c|c|}
\hline
${\rm Interval}\ I$ & $[2.0000,2.0471)$
& $[3.9801,3.9993)$ & $[4.0007,4.0952)$\\
\hline
${\rm card}\big\{{\bf z}\in\mathbb{Z}^8:\ Q_8({\bf z})\in I\big\}$ & $240$ & $270$ & $1890$\\
\hline
\end{tabular}

\vspace{0.3cm}
\centerline{Table 5. The distribution of small values
of $Q_8({\bf z})$.}
\end{table}

It is clear that the lattice $\Lambda$ corresponding to $Q_{8}({\bf z})$ is very close to $E_8$, since $G_8$ is very close to $G_{E_8}$. As a comparison, we recall (see \cite[p.123]{Co01}) that the first two nonzero shells of the positive definite quadratic form $Q_{E_8}({\bf z})$ corresponding to $E_8$ are listed in Table 6 as following.

\begin{table}[ht]
\newcommand{\AD}[1]{{\begin{tabular}{@{}c@{}}#1\end{tabular}}}
\renewcommand{\arraystretch}{1.2}
\begin{tabular}{|c|c|c|}
\hline ${\rm Norm}\ N$ & $2$ & $4$\\
\hline ${\rm card}\big\{{\bf z}\in\mathbb{Z}^8:\ Q_{E_8}({\bf z})=N\big\}$ & $240$ & $2160$\\
\hline
\end{tabular}

\vspace{0.3cm}
\centerline{Table 6. The distribution of small values of $Q_{E_8}({\bf z})$.}
\end{table}
From the nine-dimensional case onward, we omit some computational details, since the computations are similar to seven- and eight-dimensional cases.

\bigskip
\noindent
{\Large\bf 3.3. The Nine-Dimensional Case}

\medskip
Let $Q_9({\bf z})={\bf z}^{\mathsf T}G_9{\bf z}$, where $G_9=N_9/d_9,\ d_9=1000000$ and
\[
\resizebox{\textwidth}{!}{%
$N_9=\begin{pmatrix}
            5990694&3994397&999306&1996276&1994499&1995764&1992886&5573&490944\\
3994397&3999166&1002570&1997189&994558&996073&1995086&1005338&-155\\
999306&1002570&2000715&-279&-5881&-3651&3260&-1433&-505825\\
1996276&1997189&-279&2000196&439&630&992200&6772&509144\\
1994499&994558&-5881&439&2000000&998427&998427&59&492568\\
1995764&996073&-3651&630&998427&2000717&999393&1215&495538\\
1992886&1995086&3260&992200&998427&999393&2000788&980&490543\\
5573&1005338&-1433&6772&59&1215&980&2000000&-487500\\
490944&-155&-505825&509144&492568&495538&490543&-487500&2000000\\
            \end{pmatrix}.$}
\]

A routine computation shows that $G_9$ is positive definite. Using the same method as in the seven- and eight-dimensional cases, we found by computer that
$$\min_{{\bf z}\in \mathbb{Z}^9\setminus\{{\bf o}\}}Q_9({\bf z})=2.$$
In other words, the shortest nonzero vectors of the corresponding
lattice $\Lambda$ have length $\sqrt2$ and
$$r(\Lambda)=1/\sqrt2.$$
Furthermore, enumeration by computer shows that
$${\rm card}\big\{{\bf z}\in \mathbb{Z}^9:\ 2\leq Q_9({\bf z})<4\big\}=2^{10}-2=1022.$$

Using Remark 4 and the algorithm described in \cite{Av01}, we computed all vertices of the Voronoi cell $V(\Lambda)$ of the lattice $\Lambda$ corresponding to $Q_9({\bf z})$ by computer. The computation shows that $V(\Lambda)$ has $1022$ facets and $2828564$ vertices, and that the maximum distance from the origin ${\bf o}$ to a vertex is
\[
\resizebox{0.97\linewidth}{!}{%
$\displaystyle
\begin{aligned}
 R(\Lambda)&=\left(\frac{5206544919914572871728660066380044426075210316292066818456997}
{4172178730847581341450672190977943275666472067315868196000000}\right)^{\frac{1}{2}}\\[0.6em]
&=1.117103358361439607447518\cdots,
\end{aligned}
$%
}
\]
which shows that
$$\gamma^*(B^9)\leq R(\Lambda)/r(\Lambda)=\sqrt2R(\Lambda)=1.57982271\cdots<\sqrt{5/2}.$$

Enumeration by computer shows that all vectors ${\bf z}\in\mathbb{Z}^9$ satisfying $2\leq Q_9({\bf z})<4.1$ fall into the four intervals listed in Table 7.

\begin{table}[ht]
\newcommand{\AD}[1]{{\begin{tabular}{@{}c@{}}#1\end{tabular}}}
\renewcommand{\arraystretch}{1.2}
\begin{tabular}{|c|c|c|c|c|}
\hline ${\rm Interval}\ I$ & $[2.0000,2.0499)$ & $[2.9782,3.0400)$ & $[3.9642,3.99992)$ & $[4.0001,4.0982)$\\
\hline ${\rm card}\big\{{\bf z}\in\mathbb{Z}^9:\ Q_9({\bf z})\in I\big\}$ & $272$ & $256$ & $494$ & $2564$\\
\hline
\end{tabular}

\vspace{0.3cm}
\centerline{Table 7. The distribution of small values of $Q_9({\bf z})$.}
\end{table}
By the definition of the laminated lattice $\Lambda_9$ (see \cite[p.157]{Co01}), a Gram matrix $G_{{\scriptscriptstyle\sqrt{1/2}}\Lambda_9}$ of the lattice $\sqrt{1/2}\Lambda_9$ is
\[
G_{{\scriptscriptstyle\sqrt{1/2}}\Lambda_{9}}=\frac{1}{2}\begin{pmatrix}
12&8&2&4&4&4&4&0&1\\
8&8&2&4&2&2&4&2&0\\
2&2&4&0&0&0&0&0&-1\\
4&4&0&4&0&0&2&0&1\\
4&2&0&0&4&2&2&0&1\\
4&2&0&0&2&4&2&0&1\\
4&4&0&2&2&2&4&0&1\\
0&2&0&0&0&0&0&4&-1\\
1&0&-1&1&1&1&1&-1&4\\
\end{pmatrix}.
\]

It is clear that the lattice $\Lambda$ corresponding to $Q_{9}({\bf z})$ is very close to $\sqrt{1/2}\Lambda_9$, since $G_9$ is very close to $G_{{\scriptscriptstyle\sqrt{1/2}}\Lambda_9}$. As a comparison, we recall (see \cite{Sl01}) that the first three nonzero shells of the positive definite quadratic form $Q_{{\scriptscriptstyle\sqrt{1/2}}\Lambda_9}({\bf z})$ corresponding to $\sqrt{1/2}\Lambda_9$ are listed in Table 8 as following.

\begin{table}[ht]
\newcommand{\AD}[1]{{\begin{tabular}{@{}c@{}}#1\end{tabular}}}
\renewcommand{\arraystretch}{1.2}
\begin{tabular}{|c|c|c|c|}
\hline ${\rm Norm}\ N$ & $2$ & $3$ & $4$\\
\hline ${\rm card}\big\{{\bf z}\in\mathbb{Z}^9:\ Q_{{\scriptscriptstyle\sqrt{1/2}}\Lambda_9}({\bf z})=N\big\}$ & $272$ & $256$ & $3058$\\
\hline
\end{tabular}

\vspace{0.3cm}
\centerline{Table 8. The distribution of small values of $Q_{{\scriptscriptstyle\sqrt{1/2}}\Lambda_9}({\bf z})$.}
\end{table}

\bigskip\noindent
{\Large\bf 3.4. The Ten-Dimensional Case}

\medskip
Let $Q_{10}({\bf z})={\bf z}^{\mathsf T}G_{10}{\bf z}$, where $G_{10}=N_{10}/d_{10},\ d_{10}=420000$ and
\[
\resizebox{\textwidth}{!}{%
$N_{10}=\begin{pmatrix}
            1120004&-560002&-560002&280001&280001&280001&-280001&-280001&-560002&-560002\\
-560002&1120004&280001&-560002&-560002&-560002&-280001&-280001&280001&280001\\
-560002&280001&1120004&-560002&-560002&-560002&0&0&560002&0\\
280001&-560002&-560002&1120004&280001&280001&0&0&-280001&0\\
280001&-560002&-560002&280001&1120004&280001&0&0&-560002&-280001\\
280001&-560002&-560002&280001&280001&1120004&560002&560002&0&0\\
-280001&-280001&0&0&0&560002&840024&419973&139984&139981\\
-280001&-280001&0&0&0&560002&419973&840036&420036&420009\\
-560002&280001&560002&-280001&-560002&0&139984&420036&840042&279998\\
-560002&280001&0&0&-280001&0&139981&420009&279998&840018\\
            \end{pmatrix}.$}
\]

A routine computation shows that $G_{10}$ is positive definite. Using the same method as in the seven- and eight-dimensional cases, we found by computer that
$$\min_{{\bf z}\in \mathbb{Z}^{10}\setminus\{{\bf o}\}}Q_{10}({\bf z})=2.$$
In other words, the shortest nonzero vectors of the corresponding
lattice $\Lambda$ have length $\sqrt2$ and
$$r(\Lambda)=1/\sqrt2.$$
Furthermore, enumeration by computer shows that
$${\rm card}\big\{{\bf z}\in \mathbb{Z}^{10}:\ 2\leq Q_{10}({\bf z})<4\big\}=2^{11}-2=2046.$$

Using Remark 4 and the algorithm described in \cite{Av01}, we computed all vertices of the Voronoi cell $V(\Lambda)$ of the lattice $\Lambda$ corresponding to $Q_{10}({\bf z})$ by computer. The computation shows that $V(\Lambda)$ has $2046$ facets and $34941906$ vertices, and that the maximum distance from the origin ${\bf o}$ to a vertex is
\[
\resizebox{0.97\linewidth}{!}{%
$\displaystyle R(\Lambda)=\Huge\left(\frac{3155918456692576369279497300613}
{2582099165255662511773639086000}\Huge\right)^{\frac{1}{2}}=1.10554500825829207768\cdots,$}
\]
which shows that
$$\gamma^*(B^{10})\leq R(\Lambda)/r(\Lambda)=\sqrt2R(\Lambda)=1.56347674\cdots<\sqrt{8/3}.$$

Enumeration by computer shows that all the $Q_{10}({\bf z})$ satisfying $2\leq Q_{10}({\bf z})<4$ fall into the four intervals listed in Table 9.

\begin{table}[ht]
\newcommand{\AD}[1]{{\begin{tabular}{@{}c@{}}#1\end{tabular}}}
\renewcommand{\arraystretch}{1.2}
\begin{tabular}{|c|c|c|c|c|c|}
\hline ${\rm Interval}\ I$ & $[2.0,2.1)$ & $[2.6,2.7)$ & $[3.3,3.4)$ & $[3.9,4.0)$\\
\hline ${\rm card}\big\{{\bf z}\in\mathbb{Z}^{10}:\ Q_{10}({\bf z})\in I\big\}$ & $240$ & $270$ & $1296$ & $240$\\
\hline
\end{tabular}

\vspace{0.3cm}
\centerline{Table 9. The distribution of small values of $Q_{10}({\bf z})$.}

\end{table}

As a comparison, we recall (see \cite{Sl02}) that the first three nonzero shells of the positive definite quadratic form $Q_{{\scriptscriptstyle\sqrt{1/2}}\Lambda_{10}}({\bf z})$ corresponding to $\sqrt{1/2}\Lambda_{10}$ are listed in Table 10 as following.

\begin{table}[ht]
\newcommand{\AD}[1]{{\begin{tabular}{@{}c@{}}#1\end{tabular}}}
\renewcommand{\arraystretch}{1.2}
\begin{tabular}{|c|c|c|c|}
\hline ${\rm Norm}\ N$ & $2$ & $3$ & $4$\\
\hline ${\rm card}\big\{{\bf z}\in\mathbb{Z}^{10}:\ Q_{{\scriptscriptstyle\sqrt{1/2}}\Lambda_{10}}({\bf z})=N\big\}$ & $336$ & $768$ & $4950$\\
\hline
\end{tabular}

\vspace{0.3cm}
\centerline{Table 10. The distribution of small values of $Q_{{\scriptscriptstyle\sqrt{1/2}}\Lambda_{10}}({\bf z})$.}
\end{table}

Combining the four examples above, we have proved the following two theorems.

\medskip
\noindent
{\bf Theorem 1.} {\it In $7$-$10$ dimensions, the upper bounds in Table 2 can be improved as
$$\gamma^*(B^7)\leq1.51582801\cdots<\sqrt{7/3},$$
$$\gamma^*(B^8)\leq1.41413038\cdots<\sqrt2,$$
$$\gamma^*(B^9)\leq1.57982271\cdots<\sqrt{5/2}$$
and
$$\gamma^*(B^{10})\leq1.56347674\cdots<\sqrt{8/3}.$$}

\medskip
\noindent
{\bf Theorem 2.} {\it For $n=7, 8, 9$ and $10$, there exists an $n$-dimensional lattice $\Lambda$ such that all classes in $\Lambda/2\Lambda$ possess representatives of norm $N<8r(\Lambda)^2$.}

\vspace{0.8cm}
\noindent
{\LARGE\bf 4. Proof of Theorem 3}

\medskip
\noindent
Let $Q_{12}({\bf z})={\bf z}^{\mathsf T}G_{12}{\bf z}$, where $G_{12}=N_{12}/4$ and

            \[
N_{12}=\begin{pmatrix}
8&0&0&-4&0&0&4&-2&-2&-2&4&-2\\
0&8&0&0&-5&0&4&-2&-2&-2&-2&4\\
0&0&8&0&0&-4&4&4&4&-2&-2&-2\\
-4&0&0&8&0&0&-2&-2&4&4&-2&-2\\
0&-5&0&0&10&0&-2&4&-2&4&-2&-2\\
0&0&-4&0&0&8&-2&-2&-2&4&4&4\\
4&4&4&-2&-2&-2&8&0&0&-4&0&0\\
-2&-2&4&-2&4&-2&0&8&0&0&-4&0\\
-2&-2&4&4&-2&-2&0&0&8&0&0&-4\\
-2&-2&-2&4&4&4&-4&0&0&8&0&0\\
4&-2&-2&-2&-2&4&0&-4&0&0&8&0\\
-2&4&-2&-2&-2&4&0&0&-4&0&0&8\\
\end{pmatrix}.
\]

A routine computation shows that $G_{12}$ is positive definite. Using the same method as in the seven- and eight-dimensional cases, we found by computer that
$$\min_{{\bf z}\in \mathbb{Z}^{12}\setminus\{{\bf o}\}}Q_{12}({\bf z})=2.$$
In other words, the shortest nonzero vectors of the corresponding
lattice \(\Lambda\) have length \(\sqrt2\) and
$$r(\Lambda)=1/\sqrt2.$$
Furthermore, enumeration by computer shows that the first five nonzero shells of $Q_{12}({\bf z})$ are the following.

\begin{table}[ht]
\newcommand{\AD}[1]{{\begin{tabular}{@{}c@{}}#1\end{tabular}}}
\renewcommand{\arraystretch}{1.2}
\begin{tabular}{|c|c|c|c|c|c|}
\hline ${\rm Norm}\ N$ & $2$ & $5/2$ & $3$ & $7/2$ & $4$\\
\hline ${\rm card}\big\{{\bf z}\in\mathbb{Z}^{12}:\ Q_{12}({\bf z})=N\big\}$ & $594$ & $320$ & $2688$ & $1728$ & $13604$\\
\hline
\end{tabular}

\vspace{0.3cm}
\centerline{Table 11. The distribution of small values of $Q_{12}({\bf z})$.}
\end{table}

\noindent
Enumeration by computer also verifies that
$$\big\{{\bf z}\ {\rm mod}\ 2:\ {\bf z}\in\mathbb{Z}^{12},\ 	Q_{12}({\bf z})\leq4\big\}=(\mathbb{Z}/2\mathbb{Z})^{12}$$
contains $4096$ vectors. In other words, every class of $\Lambda/2\Lambda$ contains a representative ${\bf v}$ satisfying the norm
$$||{\bf v}||^2\leq4.$$

As a comparison, we recall (see \cite[p.129]{Co01}) that the first three nonzero shells of the positive definite quadratic form $Q_{{\scriptscriptstyle\sqrt{1/2}}K_{12}}({\bf z})$ corresponding to $\sqrt{1/2}K_{12}$ are listed in Table 12 as following.

\begin{table}[ht]
\newcommand{\AD}[1]{{\begin{tabular}{@{}c@{}}#1\end{tabular}}}
\renewcommand{\arraystretch}{1.2}
\begin{tabular}{|c|c|c|c|}
\hline ${\rm Norm}\ N$ & $2$ & $3$ & $4$\\
\hline ${\rm card}\big\{{\bf z}\in\mathbb{Z}^{12}:\ Q_{{\scriptscriptstyle\sqrt{1/2}}K_{12}}({\bf z})=N\big\}$ & $756$ & $4032$ & $20412$\\
\hline
\end{tabular}

\vspace{0.3cm}
\centerline{Table 12. The distribution of small values of $Q_{{\scriptscriptstyle\sqrt{1/2}}K_{12}}({\bf z})$.}
\end{table}

\noindent
Therefore, we have proved the following theorem.

\medskip
\noindent
{\bf Theorem 3.} {\it There exists a $12$-dimensional lattice $\Lambda$ (not similar to $K_{12}$) such that all classes in $\Lambda/2\Lambda$ possess representatives of norm $N\leq8r(\Lambda)^2$.}

\vfill\eject
\centerline{\Large\bf Funding}

\medskip
This work is supported by the National Natural Science Foundation of China (NSFC122- 26006, NSFC11921001) and the Natural Key Research and Development Program of China (2018YFA0704701).

\bigskip
\centerline{\Large\bf Acknowledgments}

\medskip
We acknowledge that no AI is involved in this work.

\end{document}